\documentclass[12pt]{article}%
\usepackage{graphicx}
\usepackage[intlimits]{amsmath}
\usepackage{amsfonts}
\usepackage{amssymb}%
 \makeatletter
\newtheorem{theorem}{Theorem}

\newtheorem{definition}[theorem]{Definition}

\newtheorem{problem}[theorem]{Problem}
\newtheorem{proposition}[theorem]{Proposition}

\newenvironment{proof}[1][Proof]{\noindent{\textbf {#1}  }}  {\hfill$\Box$\bigskip}

\begin{document}

\title{Extremal norms of graphs and matrices}
\author{Vladimir Nikiforov \thanks{Department of Mathematical Sciences, University of
Memphis, Memphis TN 38152, USA} \thanks{Research supported by NSF Grant
DMS-0906634.}}
\maketitle

\begin{abstract}
In the recent years, the trace norm of graphs has been extensively studied
under the name of graph energy. In this paper some of this research is
extended to more general matrix norms, like the Schatten $p$-norms and the Ky
Fan $k$-norms. Whenever possible the results are given both for graphs and
general matrices. In various contexts a puzzling fact was observed: the
Schatten $p$-norms are widely different for $1\leq p<2$ and for $p\geq
2.$\bigskip


\end{abstract}

\section{Introduction}

The aim of this paper is twofold: first, to study the Schatten and the Ky Fan
norms of the adjacency matrices of graphs, and second, obtain extremal
properties of these norms starting from graph results. We shall show that some
well-known problems and results in spectral graph theory are, in fact, best
stated in terms of matrix norms. Moreover, when stated in such terms, graph
theoretical proofs may be extended to arbitrary matrices, thus providing
extremal relations about matrix norms.

Let $G$ be a graph of order $n$ and and let $\mu_{1}\left(  G\right)
\geq\cdots\geq\mu_{n}\left(  G\right)  $ be the eigenvalues of its adjacency
matrix $A\left(  G\right)  $. Gutman \cite{Gut78,Gut99} introduced the
\emph{energy} $\mathcal{E}\left(  G\right)  $ of $G$ as
\[
\mathcal{E}\left(  G\right)  =\left\vert \mu_{1}\right\vert +\cdots+\left\vert
\mu_{n}\right\vert .
\]
This parameter, initially conceived for applications in chemistry, recently
got huge attention for its own sake: more than hundred papers have been
devoted to the study $\mathcal{E}\left(  G\right)  $ and its variants. In
particular, in \cite{Nik07}, graph energy has been generalized to arbitrary
matrices, including nonsquare ones: for an $m\times n$ matrix $A,$ define the
energy $\mathcal{E}\left(  A\right)  $ of $A$ as
\begin{equation}
\mathcal{E}\left(  A\right)  =\sigma_{1}\left(  A\right)  +\cdots+\sigma
_{m}\left(  A\right)  , \label{ener}%
\end{equation}
where $\sigma_{1}\left(  A\right)  \geq\sigma_{2}\left(  A\right)  \geq\cdots$
are the singular values of $A,$ that is to say, the square roots of the
eigenvalues of $AA^{\ast},$ and $A^{\ast}$ is the Hermitian adjoint of $A.$

Since an $m\times n$ matrix $A$ has at most $\min\left(  m,n\right)  $ nonzero
singular values, $\mathcal{E}\left(  A\right)  $ is just the sum of all
singular values of $A$. Note also that the singular values of a Hermitian
matrix are the moduli of its eigenvalues, and so, if $A$ is the adjacency
matrix of a graph $G,$ then $\mathcal{E}\left(  A\right)  =\mathcal{E}\left(
G\right)  ;$ therefore, $\mathcal{E}\left(  A\right)  $ indeed generalizes the
concept of graph energy. Moreover, the definition (\ref{ener}) sheds new light
on this concept since the sum of the singular values of $A$ is the
well-studied \emph{trace }or \emph{nuclear norm }of $A.$ It seems that graph
energy is interesting precisely because the trace norm is a fundamental matrix
parameter anyway .

Although viewing $\mathcal{E}\left(  G\right)  $ as the trace norm of
$A\left(  G\right)  $ proved to be useful (see, e.g., \cite{DaSo07},
\cite{DaSo08} and \cite{SRAG10} for some applications), this direction was not
investigated much further. To outline a new line of possible research, we
first introduce some notation and definitions.

We write $\mathcal{M}_{m,n}$ for the set of complex matrices of size $m\times
n.$ Given a matrix $A=\left[  a_{ij}\right]  \in\mathcal{M}_{m,n},$ we write
$\left\vert A\right\vert _{\infty}$ for $\max_{i,j}\left\vert a_{ij}%
\right\vert ,$ and $\left\vert A\right\vert _{p}$ for $\left(  \sum
_{i,j}\left\vert a_{ij}\right\vert ^{p}\right)  ^{1/p}.$

Now, let us recall the definitions of the Schatten and the Ky Fan matrix norms:

\begin{definition}
Given $A\in\mathcal{M}_{m,n}$ and $p\geq1,$ the Schatten $p$-norm $\left\Vert
A\right\Vert _{S_{p}}$ is given by%
\[
\left\Vert A\right\Vert _{S_{p}}=\left(  \sigma_{1}^{p}\left(  A\right)
+\cdots+\sigma_{m}^{p}\left(  A\right)  \right)  ^{1/p}.
\]

\end{definition}

\begin{definition}
Given an integer $k\geq1,$ the Ky Fan $k$-norm $\left\Vert A\right\Vert
_{F_{k}}$ is given by%
\[
\left\Vert A\right\Vert _{F_{k}}=\sigma_{1}\left(  A\right)  +\cdots
+\sigma_{k}\left(  A\right)  .
\]

\end{definition}

If $G$ is a graph with adjacency matrix $A,$ for short we write $\left\Vert
G\right\Vert _{S_{p}}$, $\left\Vert G\right\Vert _{F_{k}}$ and $\left\vert
G\right\vert _{p}$ for $\left\Vert A\right\Vert _{S_{p}},$ $\left\Vert
A\right\Vert _{F_{k}}$ and $\left\vert A\right\vert _{p},$ respectively.
Clearly, since $\left\Vert G\right\Vert _{S_{1}}=\left\Vert G\right\Vert
_{F_{n}}=\mathcal{E}\left(  G\right)  ,$ both $\left\Vert G\right\Vert
_{S_{p}}$ and $\left\Vert G\right\Vert _{F_{k}}$ extend $\mathcal{E}\left(
G\right)  $. Throughout this paper we adopt $\left\Vert G\right\Vert _{S_{1}}$
and $\left\Vert A\right\Vert _{S_{1}}$ for $\mathcal{E}\left(  G\right)  $ and
$\mathcal{E}\left(  A\right)  .$

Note that, in addition to graph energy, for positive integers $k,$ the norms
$\left\Vert G\right\Vert _{S_{2k}}$ are in wide albeit implicit use in graph
theory. Indeed, using the well-known fact about any matrix $A$%
\[
\left\Vert A\right\Vert _{S_{2}}^{2}=\sigma_{1}^{2}\left(  A\right)
+\cdots+\sigma_{m}^{2}\left(  A\right)  =tr\left(  AA^{\ast}\right)
=\left\vert A\right\vert _{2}^{2},
\]
for graphs we find that
\[
\left\Vert G\right\Vert _{S_{2}}=\sqrt{2e\left(  G\right)  },
\]
where, as usual, $e\left(  G\right)  $ stands for the number of edges of $G$.
More generally, if $k$ is a positive integer and $G$ is a graph with adjacency
matrix $A,$ then%
\[
\left\Vert G\right\Vert _{S_{2k}}=\left(  tr\left(  A^{2k}\right)  \right)
^{1/2k},
\]
and so, $\left\Vert G\right\Vert _{S_{2k}}^{2k}$ is just the number of closed
walks of length $2k$ in $G$ - a well studied graph parameter.

Since particular Schatten and Ky Fan norms are already used implicitly in
graph theory, we suggest an explicit and focused study of $\left\Vert
G\right\Vert _{S_{p}}$ and $\left\Vert G\right\Vert _{F_{k}}$ for arbitrary
$p$ and $k$. In particular, the following general problem seems interesting:

\begin{problem}
Study the extrema of $\left\Vert G\right\Vert _{S_{p}}$ and $\left\Vert
G\right\Vert _{F_{k}},$ and their relations to the structure of $G$.
\end{problem}

Following this research line, below we extend some basic results about graph
energy and other related topics.

On the other hand, many sound results in spectral graph theory can be readily
extended to matrices, sometimes even to nonsquare ones. This fact prompts
another line of investigation:

\begin{problem}
Adopting techniques from graph theory, study extremal properties of
$\left\Vert A\right\Vert _{S_{p}}$ and $\left\Vert A\right\Vert _{F_{k}},$ and
their relation to the structure of $A$ when $A$ belongs to a given class of matrices.
\end{problem}

We\ give a few such results in Section \ref{sec 2}, but they are just the tip
of the iceberg. One general finding is that, in various contexts, Schatten
$p$-norms are widely different for $1\leq p<2$ and for $p\geq2;$ we have no
explanation of this fact.

The rest of the paper is organized as follows: in Section \ref{sec1} we
discuss extremal results about the Schatten norms and Ky Fan norms of graphs:
We first find the Schatten $p$-norms of almost all graphs, which turn out to
be highly concentrated; then we extend McClelland's inequality \cite{McL71}
and some results of Koolen and Moulton \cite{KoMo01}. We also extend a lower
bound due to Caporossi, Cvetkovi\'{c}, Gutman and Hansen \cite{CCGH99}, and a
recent result of Mohar \cite{Moh09} on the sum of the largest eigenvalues of a
graph. In Section \ref{sec 2} we extend some of the results in Section
\ref{sec1} to matrices as general as possible. At the end we outline some open problems.

We do not give detailed proofs of all stated results, but we will give enough
clues to complete the proofs or to find them elsewhere.

\section{\label{sec1}Extremal results about Schatten and Ky Fan norms of
graphs}

For undefined graph theoretic and matrix notation we refer the reader to
\cite{Bol98} and \cite{HoJo88}.

Our first topic is intended to provide some intuition what are the Schatten
$p$-norms of the \textquotedblleft average\textquotedblright\ graph.

\subsection{\label{sec1.1}The Schatten norms of almost all graphs}

The purpose here is to find $\left\Vert G\right\Vert _{S_{p}}$ for all
$p\geq1$ and almost all graphs $G$. We shall prove the following theorem.

\begin{theorem}
If $G$ is a a graph of order $n$, then with probability tending to $1,$%
\[
\left\Vert G\right\Vert _{S_{p}}=\left\{
\begin{array}
[c]{ll}%
\left(  \frac{1}{\sqrt{\pi}}\cdot\frac{\Gamma\left(  p/2+1/2\right)  }%
{\Gamma\left(  p/2+2\right)  }+o\left(  1\right)  \right)  ^{1/p}n^{1/p+1/2} &
\text{if }1\leq p<2;\\
\left(  1/\sqrt{2}+o\left(  1\right)  \right)  n & \text{if }p=2;\\
\left(  1/2+o\left(  1\right)  \right)  n & \text{if }p>2.
\end{array}
\right.
\]

\end{theorem}

\begin{proof}
It is clear that the Schatten $p$-norms of almost graphs can be reduced to
finding the Schatten $p$-norms of the random graph $G\left(  n,1/2\right)  ,$
which turn out to be highly concentrated. Thus, for our calculations we will
need some results about random matrices. For short we shall write a.s. instead
of \textquotedblleft with probability tending to $1$\textquotedblright. Recall
that the adjacency matrix $A\left(  n\right)  $ of the random graph $G\left(
n,1/2\right)  $ is a random symmetric matrix with zero diagonal, whose entries
$a_{ij}$ are independent random variables with mean $\mathbb{E}\left(
a_{ij}\right)  =1/2,$ variance $Var\left(  a_{ij}^{2}\right)  =1/4,$ and
$\mathbb{E}\left(  a_{ij}^{2k}\right)  =1/4^{k}$ for all $1\leq i<j\leq n,$
$k\geq1.$ The result of F\"{u}redi and Koml\'{o}s \cite{FuKo81} implies that,
a.s.
\begin{align}
\sigma_{1}\left(  A\left(  n\right)  \right)   &  =\left(  1/2+o\left(
1\right)  \right)  n,\text{ }\label{FK1}\\
\sigma_{2}\left(  A\left(  n\right)  \right)   &  \leq\left(  1+o\left(
1\right)  \right)  n^{1/2}. \label{FK2}%
\end{align}

In finding $\left\Vert A\left(  n\right)  \right\Vert _{S_{p}},$ we shall
distinguish three cases: $p>2,$ $p=2,$ and $1\leq p<2.$ Let first $p>2.$ For
every graph $G$ of order $n$ we have%
\[
\sigma_{1}^{p}\left(  G\right)  \leq\left\Vert G\right\Vert _{S_{p}}^{p}%
\leq\sigma_{1}^{p}\left(  G\right)  +n\sigma_{2}^{p}\left(  G\right)  .
\]
By (\ref{FK1}) and (\ref{FK2}), we see that a.s.
\[
\left(  2^{-p}+o\left(  1\right)  \right)  n^{p}\leq\left\Vert A\left(
n\right)  \right\Vert _{S_{p}}^{p}\leq\left(  2^{-p}+o\left(  1\right)
\right)  n^{p}+O\left(  n^{1+p/2}\right)  ,
\]
Thus, we see that, a.s.,%
\[
\left\Vert A\left(  n\right)  \right\Vert _{S_{p}}=\left(  1/2+o\left(
1\right)  \right)  n.
\]

When $1\leq p<2$ we apply Wigner's semicircle law \cite{Wig58} in the form
given by Arnold \cite{Arn67}, p. 263. Given the conditions on $A\left(
n\right)  ,$ and the fact that $1\leq p<2,$ one can show that a. s.
\begin{align*}
\left\Vert A\left(  n\right)  \right\Vert _{S_{p}}^{p}-\sigma_{1}^{p}\left(
A\left(  n\right)  \right)   &  =n^{p/2+1}\left(  \frac{2}{\pi}\int_{-1}%
^{1}\left\vert x\right\vert ^{p}\sqrt{1-x^{2}}dx+o\left(  1\right)  \right) \\
&  =\frac{4n}{\pi}\left(  \int_{0}^{1}x^{p}\sqrt{1-x^{2}}dx+o\left(  1\right)
\right)  .
\end{align*}
Here we use that the order of $\left\Vert A\left(  n\right)  \right\Vert
_{S_{p}}^{p}$ is $n^{p/2+1}$ and all singular values of $A\left(  n\right)  $
except the largest one are bounded by (\ref{FK2}); thus, the contribution of
$\sigma_{1}\left(  A\left(  n\right)  \right)  $ is negligible and so is the
contribution of a vanishing proportion of any other singular values of $A.$

Using \cite{HDw57}, p. 196, we see that
\[
\int_{0}^{1}x^{p}\sqrt{1-x^{2}}dx=\frac{\sqrt{\pi}}{4}\cdot\frac{\Gamma\left(
p/2+1/2\right)  }{\Gamma\left(  p/2+2\right)  },
\]
and therefore, a.s., if $1\leq p<2,$ we have
\[
\left\Vert A\left(  n\right)  \right\Vert _{S_{p}}=\left(  \frac{1}{\sqrt{\pi
}}\cdot\frac{\Gamma\left(  p/2+1/2\right)  }{\Gamma\left(  p/2+2\right)
}+o\left(  1\right)  \right)  ^{1/p}n^{1/p+1/2}.
\]

For $p=2,$ it is a well-known fact that a.s. the random graph $G\left(
n,1/2\right)  $ has $\left(  1/4+o\left(  1\right)  \right)  n^{2}$ edges and
so, a.s.,
\[
\left\Vert A\left(  n\right)  \right\Vert _{S_{2}}^{2}=\left(  1/2+o\left(
1\right)  \right)  n^{2}.
\]
This completes the proof of the theorem.
\end{proof}

In particular, for $p=1,$ we see that the energy of almost all graphs of order
$n$ is a.s.
\begin{equation}
\left(  \frac{4}{3\pi}+o\left(  1\right)  \right)  n^{3/2}, \label{aag}%
\end{equation}
We established this basic fact in \cite{Nik07}.

Recently Du et al. \cite{DLL09,DLL10} calculated the energy of the random
graph $G\left(  n,p\right)  $ for $p\neq1/2.$ However, this exercise does not
generalize the above result, as these authors claim. In fact, such claim makes
no sense since graph energy is concentrated around (\ref{aag}) and cannot take
another value.

Next we investigate upper and lower bounds on the Schatten $p$-norms for any graph.

\subsection{Extremal Schatten norms of graphs}

The material in this section is by no means exhaustive: our goal is to extend
just a few of the most important bounds on graph energy.

Let $G$ is a graph of order $n$ and size $m.$ A simple, yet rather efficient
upper bound on graph energy gives the MClelland inequality \cite{McL71}
\begin{equation}
\left\Vert G\right\Vert _{S_{1}}\leq\sqrt{2mn}. \label{MCb}%
\end{equation}
Using the Cauchy-Schwarz inequality, we easily generalize this bound for
$\left\Vert G\right\Vert _{S_{p}}$ when $p>1.$ Note, however, that the
generalizations are different for $p\leq2$ and $p>2.$

\begin{proposition}
\label{McC}If $G$ is a graph of order $n$ and size $m,$ then
\[
\left\Vert G\right\Vert _{S_{p}}\leq n^{1-p/2}\left(  2m\right)  ^{p/2}%
\]
for $1\leq p\leq2,$ and
\[
\left\Vert G\right\Vert _{S_{p}}\geq n^{1-p/2}\left(  2m\right)  ^{p/2}%
\]
for $p>2.$
\end{proposition}

As we will see in Theorem \ref{MCgen}, for arbitrary matrices these bounds are
as good as one can get, but for graphs there are improvements, some of which
are discussed below.

\subsubsection{The maximal Schatten $p$-norms for $1\leq p<2$}

Koolen and Moulton (\cite{KoMo01}) proposed more sophisticated upper bounds on
$\left\Vert G\right\Vert _{S_{1}},$ starting with the following one:

\emph{If }$G$\emph{ is a graph of order }$n,$\emph{ size }$m,$\emph{ and
spectral radius }$\mu,$\emph{ then}
\begin{equation}
\left\Vert G\right\Vert _{S_{1}}\leq\mu+\left(  n-1\right)  ^{1/2}\left(
2m-\mu^{2}\right)  ^{1/2}. \label{KMg}%
\end{equation}

The straightforward proof of Koolen and Moulton in fact gives the following generalization:

\begin{proposition}
If $G$ is a graph of order $n$ and size $m$ and spectral radius $\mu,$ then
\[
\left\Vert G\right\Vert _{S_{p}}^{p}\leq\mu^{p}+\left(  n-1\right)
^{1-p/2}\left(  2m-\mu^{2}\right)  ^{p/2}%
\]
for $1\leq p\leq2.$
\end{proposition}

Next, using that $\mu\geq2m/n$, we obtain another upper bound.

\begin{theorem}
If $1\leq p\leq2$ and $G$ is a graph of order $n$ and size $m\geq n/2,$ then
\begin{equation}
\left\Vert G\right\Vert _{S_{p}}^{p}\leq\left(  2m/n\right)  ^{p}+\left(
n-1\right)  ^{1-p/2}\left(  2m-\left(  2m/n\right)  ^{2}\right)  ^{p/2}.
\label{bo2}%
\end{equation}
Equality holds if and only if $G$ is either $\left(  n/2\right)  K_{2},$ or
$K_{n},$ or a noncomplete strongly regular graph whose nontrivial eigenvalues
are
\[
\pm\left(  \frac{2m-\left(  2m/n\right)  ^{2}}{n-1}\right)  ^{1/2}.
\]

\end{theorem}

Now, we immediately see that
\[
\left(  2m/n\right)  ^{p}+\left(  n-1\right)  ^{1-p/2}\left(  2m-\left(
2m/n\right)  ^{2}\right)  ^{p/2}<n^{p}+n^{1-p/2}n^{p}2^{-p};
\]
hence the following assertion.

\begin{theorem}
If $1\leq p<2$ and $G$ is a graph of order $n,$ then
\[
\left\Vert G\right\Vert _{S_{p}}^{p}\leq2^{-p}n^{1+p/2}+n^{p}.
\]
This bound is tight: for all $n$ there is a graph $G$ of order $n$ satisfying%
\[
\left\Vert G\right\Vert _{S_{p}}^{p}=\left(  2^{-p}+o\left(  1\right)
\right)  n^{1+p/2}.
\]

\end{theorem}

To see the validity of the last statement, consider Paley graphs for prime
$n=1\left(  \operatorname{mod}4\right)  $. For any other $n$ we add isolated
vertices to the Paley graph of largest order less than $n.$ Details of such
calculations can be found in \cite{Nik07}.

By a more careful optimization Koolen and Moulton obtained the following
precise absolute bound on the energy of an $n$-vertex graph: \emph{For every
graph }$G$ \emph{of order }$n,$
\begin{equation}
\left\Vert G\right\Vert _{S_{1}}\leq\frac{n\left(  1+\sqrt{n}\right)  }{2}.
\label{KoMo}%
\end{equation}
\emph{Equality is attained for strongly regular graphs with parameters}%
\[
\left(  n,\left(  n+\sqrt{n}\right)  /2,\left(  n+2\sqrt{n}\right)  /4,\left(
n+2\sqrt{n}\right)  /4\right)  .
\]
For further details of these strongly regular graphs see Haemers and Xiang
\cite{HaXi09}. Finding the best approximation of $\left\Vert G\right\Vert
_{S_{1}}$ for all $n$ seems a challenging open problem; a partial solution is
given in \cite{Nik07a}.

\subsubsection{The maximal Schatten $p$-norms for $p\geq2$}

Let $p\geq2$ and $G$ be a graph of order $n$ and size $m$, that is to say
\begin{equation}
\sigma_{1}^{2}\left(  G\right)  +\cdots+\sigma_{n}^{2}\left(  G\right)  =2m.
\label{rest1}%
\end{equation}
To maximize $\left\Vert G\right\Vert _{S_{p}}^{p}$ subject to (\ref{rest1}),
recall that the function $x^{p/2}$ is concave for $p\geq2,$ and so,
$\left\Vert G\right\Vert _{S_{p}}^{p}$ is maximal when $\sigma_{1}^{2}\left(
G\right)  =2m$ and $\sigma_{2}\left(  G\right)  =\cdots=\sigma_{n}\left(
G\right)  =0$. This cannot really happen for graphs, but gives acceptable asymptotics.

\begin{proposition}
If $p\geq2$ and $G$ is a graph of order $n$ and size $m,$ then
\[
\left\Vert G\right\Vert _{S_{p}}^{p}<\left(  2m\right)  ^{1/2}<n.
\]
For every $m\leq n\left(  n-1\right)  /2,$ there is a graph with $m$ edges
satisfying
\[
\left\Vert G\right\Vert _{S_{p}}^{p}>\left(  2m-O\left(  m^{1/2}\right)
\right)  ^{1/2}.
\]

\end{proposition}

To see the validity of the last statement, consider the maximal complete graph
$K_{s}$ that can be formed with at most $m$ edges. We have $s=\left\lfloor
\left(  2m\right)  ^{1/2}\right\rfloor $ and
\[
\left\Vert K_{s}\right\Vert _{S_{p}}^{p}\geq\left(  s-1\right)  ^{p}+s=\left(
2m-O\left(  m^{1/2}\right)  \right)  ^{1/2}.
\]

\subsubsection{A lower bound on Schatten $p$-norms for $p\geq1$}

In this section we obtain a lower bound on $\left\Vert G\right\Vert _{S_{p}}$
which extends the bound of Caporossi, Cvetkovi\'{c}, Gutman and Hansen
\cite{CCGH99}
\begin{equation}
\mathcal{E}\left(  G\right)  \geq2\mu_{1}\left(  G\right)  , \label{Cap}%
\end{equation}
where equality holds if an only if $G$ is a complete multipartite graph with
possibly some isolated vertices. Here the extension to $\left\Vert
G\right\Vert _{S_{p}}$ sheds new light on the original bound (\ref{Cap}) by
revealing the role of the chromatic number in this relation; obviously, this
role is concealed in (\ref{Cap}).

Indeed, recall that if $G$ is a graph with chromatic number $\chi,$ Hoffman's
inequality \cite{Hof70} gives
\begin{equation}
\left\vert \mu_{n}\left(  G\right)  \right\vert +\cdots+\left\vert \mu
_{n-\chi+2}\left(  G\right)  \right\vert \geq\mu_{1}\left(  G\right)  .
\label{Hofin}%
\end{equation}
Now we can deduce the following theorem.

\begin{theorem}
\label{Schr}If $p\geq1$ and $G$ is a graph with chromatic number $\chi,$ then
\[
\left\Vert G\right\Vert _{S_{p}}\geq\sigma_{1}\left(  1+\left(  \chi-1\right)
^{1-p}\right)  ^{1/p}.
\]
For $p=1$ equality holds if an only if $G$ is a complete $\chi$-partite graph
with possibly some isolated vertices. For $p>1$ equality holds if an only if
$G$ is a regular complete $\chi$-partite graph with possibly some isolated vertices.
\end{theorem}

We switch now to the study of the Ky Fan $k$-norms.

\subsection{The maximal Ky Fan norms for graphs}

In this subsection we study the asymptotics of the maximal Ky Fan $k$-norms
for graphs of large order. To approach the problem, let us define the
functions $\tau_{k}\left(  n\right)  $ and $\xi_{k}\left(  n\right)  $ as
\begin{align*}
\tau_{k}\left(  n\right)   &  =\max_{v\left(  G\right)  =n}\mu_{1}\left(
G\right)  +\cdots+\mu_{k}\left(  G\right) \\
\xi_{k}\left(  n\right)   &  =\max_{v\left(  G\right)  =n}\left\Vert
G\right\Vert _{F_{k}}.
\end{align*}

For large $n$ the function $\tau_{k}\left(  n\right)  $ is pretty stable as
implied by the main result in \cite{Nik06}:

\begin{theorem}
For every fixed positive integer $k,$ the limit $\tau_{k}=\lim
\limits_{n\rightarrow\infty}\tau_{k}\left(  n\right)  /n$ exists.
\end{theorem}

Following the approach of \cite{Nik06},\ one can prove an analogous assertion
for singular values:

\begin{theorem}
For every fixed positive integer $k,$ the limit $\xi_{k}=\lim
\limits_{n\rightarrow\infty}\xi_{k}\left(  n\right)  /n$ exists.
\end{theorem}

Finding $\tau_{k}\left(  n\right)  $ and $\xi_{k}\left(  n\right)  $ is not
easy for any $k\geq2,$ and even finding the limits $\tau_{k}$ and $\xi_{k}$ is
challenging. Indeed, even the simplest case $\xi_{2}$ is not known yet,
despite intensive research; here is the story: Gregory, Hershkowitz and
Kirkland \cite{GHK01} asked what is the maximal value of the spread of a graph
of order $n,$ that is to say, what is
\[
\max_{v\left(  G\right)  =n}\mu_{1}\left(  G\right)  -\mu_{n}\left(  G\right)
.
\]
This problem is still open, even asymptotically; we will not solve it here,
but we will show that it is equivalent to finding%
\[
\max_{v\left(  G\right)  =n}\left\Vert G\right\Vert _{F_{2}}.
\]
Indeed, we have
\[
\left\Vert G\right\Vert _{F_{2}}=\max\left\{  \left\vert \mu_{1}\left(
G\right)  \right\vert +\left\vert \mu_{2}\left(  G\right)  \right\vert
,\left\vert \mu_{1}\left(  G\right)  \right\vert +\left\vert \mu_{n}\left(
G\right)  \right\vert \right\}  ,
\]
and from \cite{EMNA08} it is known that every graph $G$ of order $n$
satisfies
\[
\left\vert \mu_{1}\left(  G\right)  \right\vert +\left\vert \mu_{2}\left(
G\right)  \right\vert \leq\left(  1/2+\sqrt{5/12}\right)  n<1.146n.
\]
On the other hand, in \cite{GHK01}, for every $n\geq2,$ a graph $G$ of order
$n$ is shown such that
\[
\left\vert \mu_{1}\left(  G\right)  \right\vert +\left\vert \mu_{n}\left(
G\right)  \right\vert \geq\left(  2n-1\right)  /\sqrt{3};
\]
hence, $\left\vert \mu_{1}\left(  G\right)  \right\vert +\left\vert \mu
_{n}\left(  G\right)  \right\vert >1.154n$ for sufficiently large $n.$
Therefore, we see that for $n$ large,
\[
\max_{v\left(  G\right)  =n}\mu_{1}\left(  G\right)  -\mu_{n}\left(  G\right)
=\max_{v\left(  G\right)  =n}\left\Vert G\right\Vert _{F_{2}}=\xi_{2}\left(
n\right)  .
\]

\subsubsection{The asymptotics of $\tau_{k}\left(  n\right)  $ and $\xi
_{k}\left(  n\right)  $}

Surprisingly, finding $\tau_{k}$ and $\xi_{k}$ seems a tad easier for large
$k$. Indeed, Mohar \cite{Moh09} proved the following bounds%
\[
\frac{1}{2}\left(  \frac{1}{2}+\sqrt{k}+o\left(  k^{-2/5}\right)  \right)
n\leq\tau_{k}\left(  n\right)  \leq\frac{1}{2}\left(  1+\sqrt{k}\right)  n.
\]
In \cite{Nik10}, in a different way, we proved the following theorem.

\begin{theorem}
\label{tMoh}If $n$ and $k$ are integers such that $n\geq k\geq2,$ then
\[
\frac{1}{2}\left(  \frac{1}{2}+\sqrt{k}+o\left(  k^{-2/5}\right)  \right)
n\leq\xi_{k}\left(  n\right)  \leq\frac{1}{2}\left(  1+\sqrt{k}\right)  n.
\]

\end{theorem}

In fact, without extra effort, one can prove a more general upper bound about
$\left(  0,1\right)  $-matrices, which we give shortly. To describe the case
of equality we need the following definition: a matrix $A$ is called
\emph{plain} if the all one vectors $\mathbf{j}_{m}\in\mathbb{R}^{m}$ and
$\mathbf{j}_{n}\in\mathbb{R}^{n}$ are singular vectors to $\sigma_{1}\left(
A\right)  ,$ that is to say, $\sigma_{1}\left(  A\right)  =\left\langle
\mathbf{j}_{m},A\mathbf{j}_{n}\right\rangle /\sqrt{mn}.$

\begin{theorem}
\label{tNik}Let $n\geq$ $m\geq k\geq1$ be integers. If $A$ is a $\left(
0,1\right)  $-matrix of size $m\times n,$ then%
\begin{equation}
\left\Vert A\right\Vert _{F_{k}}\leq\frac{1}{2}\left(  1+\sqrt{k}\right)
\sqrt{mn}. \label{in2}%
\end{equation}
Equality holds in (\ref{in1}) if and only if the matrix $J_{m,n}-2A$ is plain
and has exactly $k$ nonzero singular values which are equal. In particular, if
$m=n=k,$ equality holds if and only if $J_{m}-2A$ is a plain Hadamard matrix.
\end{theorem}

Since for every Hermitian matrix $A$ of size $n\geq k,$
\[
\mu_{1}\left(  A\right)  +\cdots+\mu_{k}\left(  A\right)  \leq\left\Vert
A\right\Vert _{F_{k}},
\]
we see that $\tau_{k}\left(  n\right)  \leq\xi_{k}\left(  n\right)  ,$ and so
the upper bound in Theorem \ref{tMoh} implies Mohar's upper bound. The
advantage of Theorem \ref{tNik} is that this bound is extended to non-square
matrices, where eigenvalues are not applicable at all; moreover, along the
same lines, Theorem \ref{mo3} below gives an analogous statement for arbitrary
nonnegative matrices. Let us also point out that Theorem \ref{tNik} neatly
proves and generalizes the result of Koolen and Moulton (\ref{KoMo}).

There are infinitely many cases when $\xi_{k}\left(  n\right)  $ attains the
upper bound of Theorem \ref{tMoh}. Indeed, let $k$ be such that there exists a
graph $G$ of order $k$ for which the Koolen and Moulton bound (\ref{KoMo}) is
an equality. Then, for every $n\geq1,$ blowing-up $G$ by a coefficient $n,$
and calculating the Ky Fan $k$-norm of the resulting graph, we obtain%
\[
\xi_{k}\left(  kn\right)  =\frac{1}{2}\left(  1+\sqrt{k}\right)  kn.
\]
In particular, it is known that for $k=4$ the complete graph $K_{4}$ satisfies%
\[
\left\Vert K_{4}\right\Vert _{F_{4}}=6=\frac{1}{2}\left(  1+\sqrt{4}\right)
4,
\]
and so, $\xi_{4}\left(  4n\right)  =6n$ and\ $\xi_{4}=3/2.$ Other known $k$ of
this type are $k=4m^{4}$ for all integer $m>1,$ (see \cite{HaXi09}), but there
are also others, like $k=16$ and $k=36.$

Clearly, following the blow-up idea, for infinitely many triples $k,m,n$ one
can construct matrices attaining equality in (\ref{in2})$.$

\subsection{The Ky Fan norms and chromatic numbers}

Using Hoffman's inequality (\ref{Hofin}), we obtain another extension of the
bound (\ref{Cap}):

\begin{theorem}
If $G$ is a graph with chromatic number $\chi,$ then%
\begin{equation}
\left\Vert G\right\Vert _{F_{\chi}}\geq2\sigma_{1}\left(  G\right)  .
\label{in1}%
\end{equation}

\end{theorem}

Unlike Theorem \ref{Schr}, it is not easy to characterize when equality holds
in (\ref{in1}). One possible case is when $G$ is a complete $\chi$-partite
graph with some isolated vertices, but there are many other constructions,
some of which are rather complicated.

\section{\label{sec 2}Some bounds on matrix norms}

Recall that an $n\times n$ matrix $A=\left[  a_{ij}\right]  $ is called
\emph{Hadamard} if $\left\vert a_{ij}\right\vert =1,$ for all $1\leq i\leq n,$
$1\leq j\leq n$ and $AA^{\ast}=nI_{n}.$ It is known that there are no real
Hadamard matrices of some orders, but there are complex Hadamard matrices of
any order, for example, the $n\times n$ discrete Fourier transform matrix
$A=\left[  a_{ij}\right]  $ defined by
\[
a_{kj}=\exp\left(  2\pi i\left(  k-1\right)  \left(  j-1\right)  /n\right)
,\text{ \ }1\leq k\leq n\text{, }1\leq j\leq n.
\]

Below we will often meet a class of matrices akin to Hadamard matrices: for
integral $n\geq m\geq2,$ let $\mathsf{Had}_{m,n}$ denote the class of $m\times
n$ matrices whose entries have the same modulus, and whose rows are pairwise
orthogonal. Clearly if $A\in\mathsf{Had}_{m,n},$ then its row vectors have the
same length, and so $\mathsf{Had}_{n,n}$ is just the set of all scalar
multiples of Hadamard matrices. Note that $\mathsf{Had}_{m,n}$ contains
complex matrices for any $m$ and $n,$ since if $A\in\mathsf{Had}_{n,n},$ then
any $m\times n$ submatrix of $A$ belongs to $\mathsf{Had}_{m,n}.$

\subsection{Generalizing McClelland's bound}

Our first goal is to generalize McClelland's bound (\ref{MCb}) to arbitrary
matrices. Using Jensen's inequality, we immediately obtain the following
general bound.

\begin{theorem}
\label{MCgen}Let $n\geq m\geq1$ and let $A\in M_{m,n}.$ If $1\leq p\leq q,$
then
\begin{equation}
m^{-1/p}\left\Vert A\right\Vert _{S_{p}}\leq m^{-1/q}\left\Vert A\right\Vert
_{S_{q}}. \label{MCm}%
\end{equation}
Equality holds in (\ref{MCm}) if and only if the rows of $A$ are pairwise
orthogonal vectors of equal length.
\end{theorem}

The last statement can be proved by using the singular value decomposition of
$A$ to show that $AA^{\ast}$ is a scalar multiple of the identity matrix
$I_{m.}$.

Taking $A$ to be the adjacency matrix of a graph and setting $q=2$, we obtain
Proposition \ref{McC}. Restricting, in addition, the range of $p$, we obtain
an absolute bound on $\left\Vert A\right\Vert _{S_{p}}.$

\begin{theorem}
\label{MCabs}Let $n\geq m\geq1$ and let $A\in M_{m,n}.$ If $1\leq p\leq2,$
then
\begin{equation}
\left\Vert A\right\Vert _{S_{p}}\leq m^{1/p}n^{1/2}\left\vert A\right\vert
_{\infty}. \label{MCma}%
\end{equation}
Equality holds in (\ref{MCma}) if and only if $A\in\mathsf{Had}_{m,n}.$
\end{theorem}

For arbitrary matrices Theorems \ref{MCgen} and \ref{MCabs} are as good as one
can get, but for specific classes of matrices they can be improved further.
For instance, if $A$ is a nonzero nonnegative matrix, (\ref{MCm}) can be
improved as shown in the following subsection.

\subsection{Generalizing the results of Koolen and Moulton}

We start with a fine tuning of (\ref{MCma}), which is valid for arbitrary matrices.

\begin{proposition}
Let $n\geq m\geq1$ and let $A\in M_{m,n}.$ If $1\leq p\leq q,$ then
\begin{equation}
\left\Vert A\right\Vert _{S_{p}}^{p}\leq\sigma_{1}^{p}\left(  A\right)
+\left(  m-1\right)  ^{1-p/q}\left(  \left\Vert A\right\Vert _{S_{q}}%
^{q}-\sigma_{1}^{q}\left(  A\right)  \right)  ^{1/q}. \label{KMm}%
\end{equation}
Equality in (\ref{KMm}) is possible if and only if $\sigma_{2}\left(
A\right)  =\cdots=\sigma_{m}\left(  A\right)  .$
\end{proposition}

It seems technically difficult to deduce nice acceptable corollaries from
(\ref{KMm}) for general $p$ and $q.$ We shall do this for $p=1$ and $q=2.$ In
this case we see that%
\begin{equation}
\left\Vert A\right\Vert _{S_{1}}\leq\sigma_{1}\left(  A\right)  +\sqrt{\left(
m-1\right)  \left(  \left\vert A\right\vert _{2}^{2}-\sigma_{1}^{2}\left(
A\right)  \right)  }. \label{n1}%
\end{equation}
Now, noting that for nonnegative matrices $\sigma_{1}\left(  A\right)
\geq\left\vert A\right\vert _{1}/\sqrt{mn},$ we get:

\emph{If }$n\geq m\geq1$\emph{ and }$A$\emph{ is an }$m\times n$\emph{
nonnegative matrix with }$\left\vert A\right\vert _{1}\geq n\left\vert
A\right\vert _{\infty},$\emph{ then }%
\begin{equation}
\left\Vert A\right\Vert _{S_{1}}\leq\frac{\left\vert A\right\vert _{1}}%
{\sqrt{mn}}+\sqrt{\left(  m-1\right)  \left(  \left\vert A\right\vert _{2}%
^{2}-\frac{\left\vert A\right\vert _{1}^{2}}{mn}\right)  }\leq\frac{\left(
m+\sqrt{m}\right)  \sqrt{n}}{2}\left\vert A\right\vert _{\infty}. \label{n2}%
\end{equation}
Inequalities (\ref{n1}) and (\ref{n2}) were given first in \cite{Nik07};
obviously they extend the corresponding inequalities of Koolen and Moulton
\cite{KoMo01} for graphs. One question that has been omitted in \cite{Nik07}
is: \emph{When do we have the equality}
\begin{equation}
\left\Vert A\right\Vert _{S_{1}}=\frac{\left(  m+\sqrt{m}\right)  \sqrt{n}}%
{2}\left\vert A\right\vert _{\infty}. \label{KMgen}%
\end{equation}
For simplicity assume that $\left\vert A\right\vert _{\infty}=1.$ Analyzing
the proof in \cite{Nik07}, one can deduce that $A$ must be a $\left(
0,1\right)  $-matrix with
\[
\sigma_{1}\left(  A\right)  =\frac{\left(  1+\sqrt{m}\right)  \sqrt{n}}%
{2}\text{, \ }\sigma_{2}\left(  A\right)  =\cdots=\sigma_{m}\left(  A\right)
=\frac{\sqrt{n}}{2}.
\]
After some calculations, we find that $\left\vert A\right\vert _{1}=\left(
m+\sqrt{m}\right)  n/2,$ and $A$ has equal column sums and equal row sums. For
some values of $m$ and $n$ it is possible to construct matrices satisfying
(\ref{KMgen}. For instance, let $B$ be the adjacency matrix of a graph of
order $m$ having maximum energy $\left(  1+\sqrt{m}\right)  m/2.$ Consider the
$m\times2m$ block matrix $A=\left(  B,B\right)  .$ After some algebra, we see
that $A$ satisfies (\ref{KMgen}).

It turns out that nonnegative matrices satisfying (\ref{KMgen}) can be
characterized in another way

\begin{theorem}
Let $n\geq m\geq1.$ A nonnegative $m\times n$ matrix $A$ satisfies
\[
\left\Vert A\right\Vert _{S_{1}}=\frac{\left(  m+\sqrt{m}\right)  \sqrt{n}}%
{2}\left\vert A\right\vert _{\infty}%
\]
if and only if the matrix $2A-J_{m,n}$ is plain and belongs to $\mathsf{Had}%
_{m,n}$.
\end{theorem}

\subsection{The maximal Ky Fan norms of matrices}

Here, the aim is to extend Theorem \ref{tMoh} to arbitrary matrices. Using the
Cauchy- Schwarz inequality, we obtain the following theorem.

\begin{theorem}
\label{mo1}Let $n\geq m\geq k\geq1.$ If $A\in M_{m,n}$ then
\[
\left\Vert A\right\Vert _{F_{k}}\leq\sqrt{k}\left\vert A\right\vert _{2}.
\]
Equality holds for infinitely many types of matrices.
\end{theorem}

To prove the last statement, let $q\geq k$ and let $B\in$ be a $k\times q$
matrix whose rows are pairwise orthogonal vectors of equal length, that is to
say, all $k$ singular values of $B$ are equal. Writing $\otimes$ for the
Kronecker product and setting $A=B\otimes J_{r,s},$ we see that $A$ has
exactly $k$ nonzero singular values, which are equal, and so $\left\Vert
A\right\Vert _{F_{k}}\leq\sqrt{k}\left\vert A\right\vert _{2}$.

\begin{theorem}
\label{mo2}Let $n\geq m\geq k\geq1.$ If $A\in M_{m,n}$ then
\[
\left\Vert A\right\Vert _{F_{k}}\leq\sqrt{kmn}\left\vert A\right\vert
_{\infty}.
\]
Equality holds for infinitely many types of matrices.
\end{theorem}

As above, to prove the last statement, let $q\geq k,$ let $B\in\mathsf{Had}%
_{k,q},$ and set $A=B\otimes J_{r,s}.$ Since $A$ has exactly $k$ nonzero
singular values, which are equal, and since all entries of $A$ are equal in
absolute value, we have $\left\Vert A\right\Vert _{F_{k}}=\sqrt{k}\left\vert
A\right\vert _{2}=\sqrt{kmn}\left\vert A\right\vert _{\infty},$ as claimed.

These bounds are as good as one can get, but for nonnegative matrices there is
a slight improvement.

\begin{theorem}
\label{mo3}Let $n\geq m\geq k\geq1.$ If $A\in M_{m,n}$ is a nonnegative matrix
then then
\[
\left\Vert A\right\Vert _{F_{k}}\leq\frac{1}{2}\left(  1+\sqrt{k}\right)
\sqrt{mn}\left\vert A\right\vert _{\infty}%
\]
Equality holds in (\ref{in1}) if and only if $A$ is a scalar multiple of a
$\left(  0,1\right)  $-matrix, the matrix $\left\vert A\right\vert _{\infty
}\left(  J_{m,n}-2A\right)  $ is plain and has exactly $k$ nonzero singular
values which are equal.
\end{theorem}

\bigskip

\textbf{Some open problems}\bigskip

(1) Find the extrema of $\left\Vert G\right\Vert _{S_{p}}$ and $\left\Vert
G\right\Vert _{F_{k}}$ when $G$ belongs to a given monotone or hereditary property.

(2) Same problem for matrices.

(3) Find the best approximation of $\max\{\left\Vert G\right\Vert _{S_{p}%
}:v\left(  G\right)  =n\}$ for all $n$ and all $1\leq p\leq2.$

(4) Find necessary and sufficient conditions for equality in Theorems
\ref{mo1}, \ref{mo2}, and \ref{mo3}.

\end{document}